\documentclass[12pt]{article}
\usepackage{amsfonts}
\usepackage{amssymb}
\usepackage{graphicx}
\usepackage{amsmath}

\input{tcilatex}

\begin{document}

\begin{center}
{\Large On Picard Value Problem of Some Difference Polynomials}

\quad

\textbf{Zinel\^{a}abidine} \textbf{LATREUCH and Benharrat BELA\"{I}DI%
\footnote{%
Corresponding author} }

\quad

\textbf{Department of Mathematics }

\textbf{Laboratory of Pure and Applied Mathematics }

\textbf{University of Mostaganem (UMAB) }

\textbf{B. P. 227 Mostaganem-(Algeria)}

\textbf{z.latreuch@gmail.com}

\textbf{benharrat.belaidi@univ-mosta.dz}

\quad
\end{center}

\noindent \textbf{Abstract. }In this paper, we study the value distribution
of zeros of certain nonlinear difference polynomials of entire functions of
finite order.

\quad



\section{Introduction and Results}

\noindent Throughout this paper, we assume that the reader is familiar with
the fundamental results and the standard notations of the Nevanlinna's value
distribution theory $($ $\left[ 10\right] ,$ $\left[ 13\right] )$. In
addition, we will use $\rho \left( f\right) $ to denote the order of growth
of $f$, we say that a meromorphic function $a\left( z\right) $ is a small
function of $f\left( z\right) $ if $T\left( r,a\right) =S\left( r,f\right) ,$
where $S\left( r,f\right) =o\left( T\left( r,f\right) \right) ,$ as $%
r\rightarrow +\infty $ outside of a possible exceptional set of finite
logarithmic measure, we use $S\left( f\right) $ to denote the family of all
small functions with respect to $f\left( z\right) $. For a meromorphic
function $f\left( z\right) ,$ we define its shift by $f_{c}\left( z\right)
=f\left( z+c\right) .$

\quad

\noindent \qquad In 1959, Hayman proved in $\left[ 11\right] $ that if $f$
is a transcendental entire function, then $f^{n}f^{\prime }$ assume every
nonzero complex number infinitely many times, provided that $n\geq 3.$
Later, Hayman $\left[ 12\right] $ conjectured that this result remains to be
valid when $n=1$ and $n=2$. Then Mues $\left[ 18\right] $ confirmed the case
when $n=2$ and Bergweiler-Eremenko $\left[ 2\right] $ and Chen-Fang $\left[ 3%
\right] $ confirmed the case when $n=1$, independently. Since then, there
are many research publications $\left( \text{see }\left[ 17\right] \right) $
regarding this type of Picard-value problem$.$ In 1997, Bergweiler obtained
the following result.

\quad

\noindent \textbf{Theorem A. }$\left( \left[ 1\right] \right) $\textbf{\ }%
\textit{If }$f$\textit{\ is a transcendental meromorphic function of finite
order and }$q$\textit{\ is a not identically zero polynomial, then }$%
ff^{\prime }-$\textit{\ }$q$\textit{\ has infinitely many zeros.}

\quad

\noindent In 2007, Laine and Yang studied the difference analogue of
Hayman's theorem and proved the following result.

\quad

\noindent \textbf{Theorem B. }$\left( \left[ 14\right] \right) $\textbf{\ }%
\textit{Let }$f\left( z\right) $\textit{\ be a transcendental entire
function of finite order, and }$c$\textit{\ be a nonzero complex constant.
Then for }$n\geq 2,$\textit{\ }$f^{n}\left( z\right) f\left( z+c\right) $%
\textit{\ assume every non-zero value }$a\in 
\mathbb{C}
$\textit{\ infinitely often.}

\quad

\noindent \qquad In the same paper, Laine and Yang showed that Theorem B
does not remain valid for the case $n=1.$ Indeed, take $f\left( z\right)
=e^{z}+1.$ Then%
\begin{equation*}
f\left( z\right) f\left( z+\pi i\right) -1=\left( 1+e^{z}\right) \left(
1-e^{z}\right) -1=-e^{2z}.
\end{equation*}%
After their, a stream of studies on the value distribution of nonlinear
difference polynomials in $f$ has been launched and many related results
have been obtained, see e.g. $\left[ 5,14,15,16\right] .$ For example, Liu
and Yang improved the previous result and obtained the following.

\quad

\noindent \textbf{Theorem C. }$\left( \left[ 15\right] \right) $\textbf{\ }%
\textit{Let }$f\left( z\right) $\textit{\ be a transcendental entire
function of finite order, and }$c$\textit{\ be a nonzero complex constant.
Then for }$n\geq 2,$\textit{\ }$f^{n}\left( z\right) f\left( z+c\right)
-p\left( z\right) $\textit{\ has infinitely many zeros, where }$p\left(
z\right) \not\equiv 0$\textit{\ is a polynomial in }$z.$

\quad

\noindent Hence, it is natural to ask: \textit{What can be said about the
value distribution of }$f\left( z\right) f\left( z+c\right) -q\left(
z\right) ,$\textit{\ when }$f$\textit{\ is a transcendental meromorphic
function and }$q$\textit{\ be a not identically zero small function of }$f$%
\textit{? \ }In this paper, as an attempt in resolving this question, we
obtain the following results.

\quad

\noindent \textbf{Theorem 1.1} \textit{Let }$f$\textit{\ be a transcendental
entire function of finite order, let }$c_{1},$\textit{\ }$c_{2}$\textit{\ be
two nonzero complex numbers such that }$f\left( z+c_{1}\right) \not\equiv
f\left( z+c_{2}\right) $\textit{\ and }$q$\textit{\ be not identically zero
polynomial}$.$\textit{\ Then }$f\left( z\right) f\left( z+c_{1}\right)
-q\left( z\right) $\textit{\ and }$f\left( z\right) f\left( z+c_{2}\right)
-q\left( z\right) $\textit{\ at least one of them has infinitely many zeros.}

\quad

\noindent The following corollary arises directly from Theorem 1.1 and
Theorem C.

\quad

\noindent \textbf{Corollary 1.1 }\textit{Let }$n\geq 1$ \textit{be an
integer and let }$c_{1},$ $c_{2}$ $\left( c_{1}c_{2}\neq 0\right) $ \textit{%
be two distinct complex numbers. Let }$\alpha ,\beta ,$ $p_{1},$ $p_{2}$ 
\textit{and }$q$ $\left( \not\equiv 0\right) $ \textit{be nonconstant
polynomials. If }$f$ \textit{is a finite order transcendental entire
solution of }%
\begin{equation*}
\left\{ 
\begin{array}{l}
f^{n}\left( z\right) f\left( z+c_{1}\right) -q\left( z\right) =p_{1}\left(
z\right) e^{\alpha \left( z\right) } \\ 
f^{n}\left( z\right) f\left( z+c_{2}\right) -q\left( z\right) =p_{2}\left(
z\right) e^{\beta \left( z\right) }%
\end{array}%
\right. ,
\end{equation*}%
\textit{then, }$n=1$\textit{\ and }$f$\textit{\ must be a periodic function
of period }$c_{1}-c_{2}.$

\section{Some lemmas}

\noindent \qquad The following lemma is an extension of the difference
analogue of the Clunie lemma obtained by Halburd and Korhonen $\left[ 8%
\right] $.

\quad

\noindent \textbf{Lemma 2.1} $\left[ 4\right] $ \textit{Let }$f\left(
z\right) $\textit{\ be a non-constant, finite order meromorphic solution of }%
\begin{equation*}
f^{n}P\left( z,f\right) =Q\left( z,f\right) ,
\end{equation*}%
\textit{where }$P\left( z,f\right) ,$\textit{\ }$Q\left( z,f\right) $\textit{%
\ are difference polynomials in }$f\left( z\right) $\textit{\ with
meromorphic coefficients }$a_{j}\left( z\right) $\textit{\ }$\left(
j=1,\cdots ,s\right) ,$\textit{\ and let }$\delta <1.$\textit{\ If the
degree of }$Q\left( z,f\right) $ \textit{as a polynomial in }$f\left(
z\right) $ \textit{and its shifts is at most }$n,$ \textit{then}%
\begin{equation*}
m\left( r,P\left( z,f\right) \right) =o\left( \frac{T\left( r+\left\vert
c\right\vert ,f\right) }{r^{\delta }}\right) +o\left( T\left( r,f\right)
\right) +O\left( \overset{s}{\underset{j=1}{\sum }}m\left( r,a_{j}\right)
\right) .
\end{equation*}%
\textit{for all }$r$\textit{\ outside an exceptional set of finite
logarithmic measure}.

\quad 

\noindent \textbf{Lemma 2.2 }$\left[ 6\right] $ \textit{Let }$f\left(
z\right) $\textit{\ be a non-constant, finite order meromorphic function and
let }$c\neq 0$ \textit{be an arbitrary complex number. Then}%
\begin{equation*}
T\left( r,f\left( z+c\right) \right) =T\left( r,f\left( z\right) \right)
+S\left( r,f\right) .
\end{equation*}%
\textbf{Lemma 2.3}\textit{\ }$\left[ 7\right] $ \textit{Let }$f\left(
z\right) $\textit{\ be a transcendental meromorphic function of finite order 
}$\rho ,$\textit{\ and let }$\varepsilon >0$\textit{\ be a given constant.
Then, there exists a set }$E_{0}\subset \left( 1,+\infty \right) $\textit{\
that has finite logarithmic measure, such that for all }$z$\textit{\
satisfying }$\left\vert z\right\vert \notin E_{0}\cup \left[ 0,1\right] ,$%
\textit{\ and for all }$k,j,$\textit{\ }$0\leq j<k,$\textit{\ we have}%
\begin{equation*}
\left\vert \frac{f^{\left( k\right) }\left( z\right) }{f^{\left( j\right)
}\left( z\right) }\right\vert \leq \left\vert z\right\vert ^{\left(
k-j\right) \left( \rho -1+\varepsilon \right) }.
\end{equation*}

\noindent \qquad The following lemma is the lemma of the logarithmic
derivative.

\quad

\noindent \textbf{Lemma 2.4 }$\left[ 10\right] $ \textit{Let }$f$ \textit{be
a meromorphic function and let }$k\in 
\mathbb{N}
.$ \textit{Then}%
\begin{equation*}
m\left( r,\frac{f^{\left( k\right) }}{f}\right) =S\left( r,f\right) ,
\end{equation*}%
\textit{where }$S\left( r,f\right) =O\left( \log T\left( r,f\right) +\log
r\right) ,$ \textit{possibly outside a set }$E_{1}\subset \left[ 0,+\infty
\right) $ \textit{of a finite linear measure}$.$ \textit{If }$f$ \textit{is
of finite order of growth, then}%
\begin{equation*}
m\left( r,\frac{f^{\left( k\right) }}{f}\right) =O\left( \log r\right) .
\end{equation*}

\noindent \qquad The following lemma is a difference analogue of the lemma
of the logarithmic derivative for finite order meromorphic functions.

\quad

\noindent \textbf{Lemma 2.5 }$\left[ 6,8,9\right] $ \textit{Let }$\eta
_{1},\eta _{2}$ \textit{be two arbitrary complex numbers such that }$\eta
_{1}\neq \eta _{2}$ \textit{and let }$f\left( z\right) $ \textit{be a finite
order meromorphic function. Let }$\sigma $ \textit{be the order of }$f\left(
z\right) $. \textit{Then for each }$\varepsilon >0,$ \textit{we have}%
\begin{equation*}
m\left( r,\frac{f\left( z+\eta _{1}\right) }{f\left( z+\eta _{2}\right) }%
\right) =O\left( r^{\sigma -1+\varepsilon }\right) .
\end{equation*}%
\textbf{Lemma 2.6}\textit{\ Let }$f\left( z\right) $\textit{\ be a
transcendental meromorphic solution of the system}%
\begin{equation}
\left\{ 
\begin{array}{l}
f\left( z\right) f\left( z+c_{1}\right) -q\left( z\right) =p_{1}\left(
z\right) e^{\alpha \left( z\right) }, \\ 
f\left( z\right) f\left( z+c_{2}\right) -q\left( z\right) =p_{2}\left(
z\right) e^{\beta \left( z\right) },%
\end{array}%
\right.  \tag{2.1}
\end{equation}%
\textit{where }$\alpha ,\beta $\textit{\ are polynomials and }$p_{1},$ $%
p_{2},$ $q$ \textit{are not identically zero rational functions. If }$%
N\left( r,f\right) =S\left( r,f\right) ,$ \textit{then }%
\begin{equation*}
\deg \alpha =\deg \beta =\deg \left( \alpha +\beta \right) =\rho \left(
f\right) >0.
\end{equation*}%
\textit{Proof}\textbf{. }First,\textbf{\ }we prove that $\deg \alpha =\rho
\left( f\right) $ and by the same we can deduce that $\deg \beta =\rho
\left( f\right) .$ It's clear from $\left( 2.1\right) $ that $\deg \alpha
\leq \rho \left( f\right) .$ Suppose that $\deg \alpha <\rho \left( f\right)
,$ this means that 
\begin{equation}
f\left( z\right) f\left( z+c_{1}\right) :=F=q\left( z\right) +p_{1}\left(
z\right) e^{\alpha \left( z\right) }\in S\left( f\right) .  \tag{2.2}
\end{equation}%
Applying Lemma 2.1 and Lemma 2.2 into $\left( 2.2\right) ,$ we obtain \ $%
T\left( r,f_{c}\right) =T\left( r,f\right) =S\left( r,f\right) $ which is a
contradiction. Assume now that $\deg \left( \alpha +\beta \right) <\rho
\left( f\right) ,$ this leads to $p_{1}p_{2}e^{\alpha +\beta }\in S\left(
f\right) .$ From this and $\left( 2.1\right) $ we have%
\begin{equation*}
f^{2}P\left( z,f\right) =p_{1}p_{2}e^{\alpha +\beta }+q^{2},
\end{equation*}%
where 
\begin{equation*}
P\left( z,f\right) =a\left( z\right) f^{2}-b\left( z\right)
\end{equation*}%
and%
\begin{equation*}
a=\frac{f_{c_{1}}}{f}\frac{f_{c_{2}}}{f},\text{ }b=q\left( \frac{f_{c_{1}}}{f%
}+\frac{f_{c_{2}}}{f}\right) .
\end{equation*}%
It's clear that $P\left( z,f\right) \not\equiv 0,$ and by using Lemma 2.1,
we get%
\begin{equation*}
m\left( r,P\left( z,f\right) \right) =S\left( r,f\right)
\end{equation*}%
which leads to 
\begin{equation*}
2T\left( r,f\right) =m\left( r,\frac{b\left( z\right) +P\left( z,f\right) }{%
a\left( z\right) }\right) =S\left( r,f\right)
\end{equation*}%
which is a contradiction. Hence, $\deg \left( \alpha +\beta \right) =\deg
\alpha =\deg \beta .$ Finally, by using Lemma 2.1, it's easy to see that
both of $\alpha $ and $\beta $ are nonconstant polynomials.

\section{Proof of Theorem 1.1}

\noindent We shall prove this theorem by contradiction. Suppose contrary to
our assertion that both of $f\left( z\right) f\left( z+c_{1}\right) -q\left(
z\right) $ and $f\left( z\right) f\left( z+c_{2}\right) -q\left( z\right) $
have finitely many zeros. Then, there exist four polynomials $\alpha $, $%
\beta ,$ $p_{1}$ and $p_{2}$ such that%
\begin{equation}
f\left( z\right) f\left( z+c_{1}\right) -q\left( z\right) =p_{1}\left(
z\right) e^{\alpha \left( z\right) }  \tag{3.1}
\end{equation}%
and 
\begin{equation}
f\left( z\right) f\left( z+c_{2}\right) -q\left( z\right) =p_{2}\left(
z\right) e^{\beta \left( z\right) }.  \tag{3.2}
\end{equation}%
By differentiating $\left( 3.1\right) $ and eliminating $e^{\alpha },$ we get%
\begin{equation}
A_{1}ff_{c_{1}}-f^{\prime }f_{c_{1}}-ff_{c_{1}}^{\prime }=B_{1},  \tag{3.3}
\end{equation}%
where $A_{1}=\frac{p_{1}^{\prime }}{p_{1}}+\alpha ^{\prime },$ $B_{1}=\left( 
\frac{p_{1}^{\prime }}{p_{1}}+\alpha ^{\prime }\right) q-q^{\prime }.$ By
Lemma 2.6 we have 
\begin{equation*}
\deg \alpha =\deg \beta =\deg \left( \alpha +\beta \right) =\rho \left(
f\right) >0.
\end{equation*}%
Now, we prove that $A_{1}\not\equiv 0.$ To show this, we suppose the
contrary. Then, there exists a constant $A$ such that $A=p_{1}\left(
z\right) e^{\alpha },$ which implies the contradiction $\deg \alpha =\rho
\left( f\right) =0.$ By the same, we can prove that $B_{1}\not\equiv 0.$ By
the same arguments as above, $\left( 3.2\right) $ gives%
\begin{equation}
A_{2}ff_{c_{2}}-f^{\prime }f_{c_{2}}-ff_{c_{2}}^{\prime }=B_{2},  \tag{3.4}
\end{equation}%
where $A_{2}=\frac{p_{2}^{\prime }}{p_{2}}+\beta ^{\prime }$ and $%
B_{2}=\left( \frac{p_{2}^{\prime }}{p_{2}}+\beta ^{\prime }\right)
q-q^{\prime }.$ Obviously, $A_{2}\not\equiv 0$ and $B_{2}\not\equiv 0.$
Dividing both sides of $\left( 3.3\right) $ and $\left( 3.4\right) $ by $%
f^{2},$ we get for each $\varepsilon >0$ 
\begin{equation*}
2m\left( r,\frac{1}{f}\right) \leq m\left( r,\frac{f_{c_{i}}}{f}\right)
+m\left( r,\frac{f^{\prime }}{f}\frac{f_{c_{i}}}{f}\right) +m\left( r,\frac{%
f_{c_{i}}^{\prime }}{f_{c_{i}}}\frac{f_{c_{i}}}{f}\right) +O\left( \log
r\right)
\end{equation*}%
\begin{equation*}
=O\left( r^{\rho -1+\varepsilon }\right) +O\left( \log r\right) =S\left(
r,f\right) .
\end{equation*}%
So, by the first fundamental theorem, we deduce that%
\begin{equation}
T\left( r,f\right) =N\left( r,\frac{1}{f}\right) +O\left( r^{\rho
-1+\varepsilon }\right) +O\left( \log r\right) .  \tag{3.5}
\end{equation}%
It's clear from $\left( 3.3\right) $ and $\left( 3.4\right) $ that any
multiple zero of $f$ is a zero of $B_{i}$ $\left( i=1,2\right) .$ Hence 
\begin{equation*}
N_{(2}\left( r,\frac{1}{f}\right) \leq N\left( r,\frac{1}{B_{i}}\right)
=O\left( \log r\right) ,
\end{equation*}%
where $N_{(2}\left( r,\frac{1}{f}\right) $ denotes the counting function of
zeros of $f$ whose multiplicities are not less than 2. It follows by this
and $\left( 3.5\right) $ that 
\begin{equation}
T\left( r,f\right) =N_{1)}\left( r,\frac{1}{f}\right) +O\left( r^{\rho
-1+\varepsilon }\right) +O\left( \log r\right) ,  \tag{3.6}
\end{equation}%
where $N_{1)}\left( r,\frac{1}{f}\right) $ is the counting function of
zeros, where only the simple zeros are considered. From $\left( 3.3\right) $
and $\left( 3.4\right) ,$ for every zero $z_{0}$ such that $f^{\prime
}\left( z_{0}\right) \neq 0$ which is not zero or pole of $B_{1}$ and $%
B_{2}, $ we have 
\begin{equation}
\left( f^{\prime }f_{c_{1}}+B_{1}\right) \left( z_{0}\right) =0  \tag{3.7}
\end{equation}%
and%
\begin{equation}
\left( f^{\prime }f_{c_{2}}+B_{2}\right) \left( z_{0}\right) =0.  \tag{3.8}
\end{equation}%
By $\left( 3.7\right) $ and $\left( 3.8\right) ,$ we obtain%
\begin{equation}
\left( B_{2}f_{c_{1}}-B_{1}f_{c_{2}}\right) \left( z_{0}\right) =0  \tag{3.9}
\end{equation}%
which means that the function $\frac{B_{2}f_{c_{1}}-B_{1}f_{c_{2}}}{f}$ has
at most a finite number of simple poles. We consider two cases:

\quad

\noindent \textbf{Case 1.} $B_{2}f_{c_{1}}-B_{1}f_{c_{2}}\not\equiv 0.$ Set%
\begin{equation}
h\left( z\right) =\frac{B_{2}f_{c_{1}}-B_{1}f_{c_{2}}}{f\left( z\right) }. 
\tag{3.10}
\end{equation}%
Then, from the lemma of logarithmic differences, we have $m\left( r,h\right)
=O\left( r^{\rho -1+\varepsilon }\right) +O\left( \log r\right) .$ On the
other hand%
\begin{equation*}
N\left( r,h\right) =N\left( r,\frac{B_{2}f_{c_{1}}-B_{1}f_{c_{2}}}{f}\right)
=N_{1)}\left( r,\frac{B_{2}f_{c_{1}}-B_{1}f_{c_{2}}}{f}\right) 
\end{equation*}%
\begin{equation*}
+O\left( r^{\rho -1+\varepsilon }\right) +O\left( \log r\right) =S\left(
r,f\right) .
\end{equation*}%
Thus, $T\left( r,h\right) =O\left( r^{\rho -1+\varepsilon }\right) +O\left(
\log r\right) =S\left( r,f\right) .$ From the equation $\left( 3.10\right) ,$
we have 
\begin{equation}
f_{c_{1}}\left( z\right) =\frac{B_{1}}{B_{2}}f_{c_{2}}\left( z\right) +\frac{%
h}{B_{2}}f\left( z\right) .  \tag{3.11}
\end{equation}%
By differentiating $\left( 3.11\right) ,$ we get%
\begin{equation}
f_{c_{1}}^{\prime }\left( z\right) =\left( \frac{h}{B_{2}}\right) ^{\prime
}f\left( z\right) +\frac{h}{B_{2}}f^{\prime }\left( z\right) +\left( \frac{%
B_{1}}{B_{2}}\right) ^{\prime }f_{c_{2}}\left( z\right) +\frac{B_{1}}{B_{2}}%
f_{c_{2}}^{\prime }\left( z\right) .  \tag{3.12}
\end{equation}%
Substituting $\left( 3.11\right) $ and $\left( 3.12\right) $ into $\left(
3.3\right) $%
\begin{equation*}
\left[ \frac{A_{1}h}{B_{2}}-\left( \frac{h}{B_{2}}\right) ^{\prime }\right]
f^{2}+\left[ -\frac{2h}{B_{2}}\right] ff^{\prime }
\end{equation*}%
\begin{equation}
+\left[ \frac{A_{1}B_{1}}{B_{2}}-\left( \frac{B_{1}}{B_{2}}\right) ^{\prime }%
\right] ff_{c_{2}}-\frac{B_{1}}{B_{2}}f^{\prime }f_{c_{2}}-\frac{B_{1}}{B_{2}%
}ff_{c_{2}}^{\prime }=B_{1}.  \tag{3.13}
\end{equation}%
Equation $\left( 3.4\right) ,$ can be rewritten as 
\begin{equation*}
-\frac{B_{1}A_{2}}{B_{2}}ff_{c_{2}}+\frac{B_{1}}{B_{2}}f^{\prime }f_{c_{2}}+%
\frac{B_{1}}{B_{2}}ff_{c_{2}}^{\prime }=-B_{1}.
\end{equation*}%
By adding this to $\left( 3.13\right) $, we get%
\begin{equation}
\left[ \frac{A_{1}h}{B_{2}}-\left( \frac{h}{B_{2}}\right) ^{\prime }\right]
f+\left[ -\frac{2h}{B_{2}}\right] f^{\prime }+\left[ \frac{A_{1}B_{1}}{B_{2}}%
-\left( \frac{B_{1}}{B_{2}}\right) ^{\prime }-\frac{B_{1}A_{2}}{B_{2}}\right]
f_{c_{2}}=0.  \tag{3.14}
\end{equation}%
Its clear that $-\frac{2h}{B_{2}}\not\equiv 0$. In order to complete the
proof of our theorem, we need to prove 
\begin{equation*}
\frac{A_{1}h}{B_{2}}-\left( \frac{h}{B_{2}}\right) ^{\prime }\not\equiv 0%
\text{ and }\frac{A_{1}B_{1}}{B_{2}}-\left( \frac{B_{1}}{B_{2}}\right)
^{\prime }-\frac{B_{1}A_{2}}{B_{2}}\not\equiv 0.
\end{equation*}%
Suppose contrary to our assertion that $\frac{A_{1}h}{B_{2}}-\left( \frac{h}{%
B_{2}}\right) ^{\prime }\equiv 0.$ Then, by the definition of $A_{1}$ and by
simple integration, we get%
\begin{equation*}
p_{1}e^{\alpha }=C_{1}\frac{h}{B_{2}},
\end{equation*}%
where $C_{1}$ is a nonzero constant. This implies that $\deg \alpha =\rho
\left( f\right) -1,$ which is a contradiction. Hence, $\frac{A_{1}h}{B_{2}}%
-\left( \frac{h}{B_{2}}\right) ^{\prime }\not\equiv 0.$ Next, we shall prove 
$\frac{A_{1}B_{1}}{B_{2}}-\left( \frac{B_{1}}{B_{2}}\right) ^{\prime }-\frac{%
B_{1}A_{2}}{B_{2}}\not\equiv 0.$ Suppose that $\frac{A_{1}B_{1}}{B_{2}}%
-\left( \frac{B_{1}}{B_{2}}\right) ^{\prime }-\frac{B_{1}A_{2}}{B_{2}}\equiv
0.$ Then we obtain%
\begin{equation*}
\frac{p_{1}}{p_{2}}e^{\alpha -\beta }=C_{2}\frac{B_{1}}{B_{2}}:=\gamma ,
\end{equation*}%
where $C_{2}$ is a nonzero constant and $\gamma $ is a small function of $f.$
From $\left( 3.1\right) $ and $\left( 3.2\right) $ we get%
\begin{equation}
f\left( f_{c_{1}}-\gamma f_{c_{2}}\right) =\left( 1-\gamma \right) q. 
\tag{3.15}
\end{equation}%
If $\gamma \not\equiv 1,$ then by applying Clunie's lemma to $\left(
3.15\right) ,$ we obtain 
\begin{equation*}
m\left( r,f_{c_{1}}-\gamma f_{c_{2}}\right) =T\left( r,f_{c_{1}}-\gamma
f_{c_{2}}\right) =S\left( r,f\right) .
\end{equation*}%
By this and $\left( 3.15\right) ,$ we have%
\begin{equation*}
T\left( r,f\right) =T\left( r,\frac{\left( 1-\gamma \right) q}{%
f_{c_{1}}-\gamma f_{c_{2}}}\right) =S\left( r,f\right) 
\end{equation*}%
which is a contradiction. If $\gamma \equiv 1,$ then we obtain the
contradiction $f_{c_{1}}\left( z\right) \equiv f_{c_{2}}\left( z\right) .$
Thus, $\frac{A_{1}B_{1}}{B_{2}}-\left( \frac{B_{1}}{B_{2}}\right) ^{\prime }-%
\frac{B_{1}A_{2}}{B_{2}}\not\equiv 0.$ From the above discussion and $\left(
3.14\right) ,$ we have%
\begin{equation}
f_{c_{2}}\left( z\right) =M\left( z\right) f\left( z\right) +N\left(
z\right) f^{\prime }\left( z\right)   \tag{3.16}
\end{equation}%
and 
\begin{equation}
f_{c_{1}}\left( z\right) =\varphi \left( z\right) f\left( z\right) +\psi
\left( z\right) f^{\prime }\left( z\right) ,  \tag{3.17}
\end{equation}%
where 
\begin{equation*}
M=\frac{\left( \frac{h}{B_{2}}\right) ^{\prime }-A_{1}\frac{h}{B_{2}}}{%
\left( A_{1}-A_{2}\right) \frac{B_{1}}{B_{2}}-\left( \frac{B_{1}}{B_{2}}%
\right) ^{\prime }},\text{ }N=\frac{\frac{2h}{B_{2}}}{\left(
A_{1}-A_{2}\right) \frac{B_{1}}{B_{2}}-\left( \frac{B_{1}}{B_{2}}\right)
^{\prime }}
\end{equation*}%
and 
\begin{equation*}
\varphi \left( z\right) =\frac{B_{1}}{B_{2}}M+\frac{h}{B_{2}},\text{ }\psi =%
\frac{B_{1}}{B_{2}}N.
\end{equation*}%
Differentiation of $\left( 3.16\right) ,$ gives%
\begin{equation}
f_{c_{2}}^{\prime }=M^{\prime }f+\left( M+N^{\prime }\right) f^{\prime
}+Nf^{\prime \prime }.  \tag{3.18}
\end{equation}%
Substituting $\left( 3.16\right) $ and $\left( 3.18\right) $ into $\left(
3.4\right) ,$ we get%
\begin{equation}
\left[ M^{\prime }-A_{2}M\right] f^{2}+\left[ N^{\prime }-A_{2}N+2M\right]
f^{\prime }f+N\left( \left( f^{\prime }\right) ^{2}+ff^{\prime \prime
}\right) =-B_{2}.  \tag{3.19}
\end{equation}%
Differentiating $\left( 3.19\right) ,$ we get%
\begin{equation*}
\left[ M^{\prime }-A_{2}M\right] ^{\prime }f^{2}+\left( 2\left[ M^{\prime
}-A_{2}M\right] +\left[ N^{\prime }-A_{2}N+2M\right] ^{\prime }\right)
f^{\prime }f
\end{equation*}%
\begin{equation}
+\left( 2N^{\prime }-A_{2}N+2M\right) \left( \left( f^{\prime }\right)
^{2}+ff^{\prime \prime }\right) +N\left( 3f^{\prime }f^{\prime \prime
}+ff^{\prime \prime \prime }\right) =-B_{2}^{\prime }.  \tag{3.20}
\end{equation}%
Suppose $z_{0}$ is a simple zero of $f$ and not a zero or pole of $B_{2}.$
Then from $\left( 3.19\right) $ and $\left( 3.20\right) ,$ we have%
\begin{equation*}
\left( Nf^{\prime }+\frac{B_{2}}{f^{\prime }}\right) \left( z_{0}\right) =0,
\end{equation*}%
\begin{equation*}
\left[ \left( 2N^{\prime }-A_{2}N+2M\right) f^{\prime }+3Nf^{\prime \prime }+%
\frac{B_{2}^{\prime }}{f^{\prime }}\right] \left( z_{0}\right) =0.
\end{equation*}%
It follows that $z_{0}$ is a zero of $\left[ B_{2}\left( 2N^{\prime
}-A_{2}N+2M\right) -B_{2}^{\prime }N\right] f^{\prime }+3B_{2}Nf^{\prime
\prime }.$ Therefore the function 
\begin{equation*}
H=\frac{\left[ 2B_{2}N^{\prime }-B_{2}A_{2}N+2B_{2}M-B_{2}^{\prime }N\right]
f^{\prime }+3B_{2}Nf^{\prime \prime }}{f}
\end{equation*}%
satisfies $T\left( r,H\right) =S\left( r,f\right) $ and%
\begin{equation}
f^{\prime \prime }=\frac{H}{3B_{2}N}f+\frac{\left[ -2B_{2}N^{\prime
}+B_{2}A_{2}N-2B_{2}M+B_{2}^{\prime }N\right] }{3B_{2}N}f^{\prime }. 
\tag{3.21}
\end{equation}%
Substituting $\left( 3.21\right) $ into $\left( 3.19\right) ,$ we get%
\begin{equation}
q_{1}f^{2}+q_{2}f^{\prime }f+q_{3}\left( f^{\prime }\right) ^{2}=-B_{2}, 
\tag{3.22}
\end{equation}%
where 
\begin{equation*}
q_{1}=M^{\prime }-A_{2}M+\frac{H}{3B_{2}},
\end{equation*}%
\begin{equation*}
q_{2}=\frac{1}{3}N^{\prime }+\frac{1}{3}\left( \frac{B_{2}^{\prime }}{B_{2}}%
-2A_{2}\right) N+\frac{4}{3}M,\text{ }q_{3}=N.
\end{equation*}%
We prove first $q_{2}\not\equiv 0.$ Suppose the contrary. Then%
\begin{equation*}
\frac{q_{2}}{q_{3}}=\frac{2}{3}\frac{N^{\prime }}{N}-\frac{1}{3}\frac{%
B_{2}^{\prime }}{B_{2}}-\frac{2}{3}\left( A_{1}+A_{2}\right) +\frac{2}{3}%
\frac{h^{\prime }}{h}=0
\end{equation*}%
which leads to%
\begin{equation*}
\alpha ^{\prime }+\beta ^{\prime }=\frac{N^{\prime }}{N}-2\frac{%
B_{2}^{\prime }}{B_{2}}+\frac{h^{\prime }}{h}-\frac{p_{1}^{\prime }}{p_{1}}-%
\frac{p_{2}^{\prime }}{p_{2}}.
\end{equation*}%
By simple integration of both sides of the above equation, we get%
\begin{equation}
p_{1}p_{2}e^{\alpha +\beta }=c\frac{N}{B_{2}^{2}}h,  \tag{3.23}
\end{equation}%
where $c$ is a nonzero constant, this leads to the contradiction $\deg
\left( \alpha +\beta \right) <\deg \alpha =\deg \beta .$ Hence, $%
q_{2}\not\equiv 0.$ Differentiating $\left( 3.22\right) ,$ we obtain 
\begin{equation}
q_{1}^{\prime }f^{2}+\left( 2q_{1}+q_{2}^{\prime }\right) f^{\prime
}f+\left( q_{2}+q_{3}^{\prime }\right) \left( f^{\prime }\right)
^{2}+q_{2}f^{\prime \prime }f+2q_{3}f^{\prime }f^{\prime \prime
}=-B_{2}^{\prime }.  \tag{3.24}
\end{equation}%
Let $z_{0}$ be a simple zero of $f$ which is not a zero or pole of $B_{2}.$
Then from $\left( 3.22\right) $ and $\left( 3.24\right) $ we have%
\begin{equation*}
\left( q_{3}f^{\prime }+\frac{B_{2}}{f^{\prime }}\right) \left( z_{0}\right)
=0,
\end{equation*}%
\begin{equation*}
\left[ \left( q_{2}+q_{3}^{\prime }\right) f^{\prime }+2q_{3}f^{\prime
\prime }+\frac{B_{2}^{\prime }}{f^{\prime }}\right] \left( z_{0}\right) =0.
\end{equation*}%
Therefore $z_{0}$ is a zero of $\left( B_{2}\left( q_{2}+q_{3}^{\prime
}\right) -B_{2}^{\prime }q_{3}\right) f^{\prime }+2B_{2}q_{3}f^{\prime
\prime }.$ Hence the function 
\begin{equation*}
R=\frac{\left( B_{2}\left( q_{2}+q_{3}^{\prime }\right) -B_{2}^{\prime
}q_{3}\right) f^{\prime }+2B_{2}q_{3}f^{\prime \prime }}{f}.
\end{equation*}%
satisfies $T\left( r,R\right) =S\left( r,f\right) $ and 
\begin{equation}
f^{\prime \prime }=\frac{R}{2B_{2}q_{3}}f+\frac{B_{2}^{\prime
}q_{3}-B_{2}\left( q_{2}+q_{3}^{\prime }\right) }{2B_{2}q_{3}}f^{\prime }. 
\tag{3.25}
\end{equation}%
Substituting $\left( 3.25\right) $ into $\left( 3.24\right) $%
\begin{equation*}
\left[ q_{1}^{\prime }+\frac{q_{2}R}{2B_{2}q_{3}}\right] f^{2}+\left[
2q_{1}+q_{2}^{\prime }+\frac{1}{2}\frac{B_{2}^{\prime }}{B_{2}}q_{2}-\frac{1%
}{2}\left( q_{2}+q_{3}^{\prime }\right) \frac{q_{2}}{q_{3}}+\frac{R}{B_{2}}%
\right] f^{\prime }f
\end{equation*}%
\begin{equation}
+\frac{B_{2}^{\prime }q_{3}}{B_{2}}\left( f^{\prime }\right)
^{2}=-B_{2}^{\prime }.  \tag{3.26}
\end{equation}%
Combining $\left( 3.26\right) $ and $\left( 3.22\right) ,$ we obtain%
\begin{equation}
\left[ q_{1}^{\prime }+\frac{q_{2}R}{2B_{2}q_{3}}-\frac{B_{2}^{\prime }}{%
B_{2}}q_{1}\right] f+\left[ 2q_{1}+q_{2}^{\prime }-\frac{1}{2}\frac{%
B_{2}^{\prime }}{B_{2}}q_{2}-\frac{1}{2}\left( q_{2}+q_{3}^{\prime }\right) 
\frac{q_{2}}{q_{3}}+\frac{R}{B_{2}}\right] f^{\prime }=0.  \tag{3.27}
\end{equation}%
From $\left( 3.27\right) ,$ we deduce that 
\begin{equation*}
q_{1}^{\prime }+\frac{q_{2}R}{2B_{2}q_{3}}-\frac{B_{2}^{\prime }}{B_{2}}%
q_{1}=0
\end{equation*}%
and 
\begin{equation*}
2q_{1}+q_{2}^{\prime }-\frac{1}{2}\frac{B_{2}^{\prime }}{B_{2}}q_{2}-\frac{1%
}{2}\left( q_{2}+q_{3}^{\prime }\right) \frac{q_{2}}{q_{3}}+\frac{R}{B_{2}}%
=0.
\end{equation*}%
By eliminating $R$ from the above two equations, we obtain%
\begin{equation}
q_{3}\left( 4q_{1}q_{3}-q_{2}^{2}\right) \frac{B_{2}^{\prime }}{B_{2}}%
+q_{2}\left( 4q_{1}q_{3}-q_{2}^{2}\right) -q_{3}\left(
4q_{1}q_{3}-q_{2}^{2}\right) ^{\prime }+q_{3}^{\prime }\left(
4q_{1}q_{3}-q_{2}^{2}\right) =0.  \tag{3.28}
\end{equation}%
Thus, equation $\left( 3.25\right) $ can be rewritten as 
\begin{equation}
f^{\prime \prime }=\left( \frac{B_{2}^{\prime }}{B_{2}}\frac{q_{1}}{q_{2}}-%
\frac{q_{1}^{\prime }}{q_{2}}\right) f+\frac{1}{2}\left( \frac{B_{2}^{\prime
}}{B_{2}}-\frac{q_{2}}{q_{3}}-\frac{N^{\prime }}{N}\right) f^{\prime }. 
\tag{3.29}
\end{equation}%
\textbf{Subcase 1.} If $4q_{1}q_{3}-q_{2}^{2}\not\equiv 0,$ then from $%
\left( 3.28\right) $ we have%
\begin{equation*}
\frac{q_{2}}{q_{3}}=\frac{\left( 4q_{1}q_{3}-q_{2}^{2}\right) ^{\prime }}{%
\left( 4q_{1}q_{3}-q_{2}^{2}\right) }-\frac{B_{2}^{\prime }}{B_{2}}-\frac{%
q_{3}^{\prime }}{q_{3}}.
\end{equation*}%
On the other hand 
\begin{equation*}
\frac{q_{2}}{q_{3}}=\frac{1}{3}\frac{N^{\prime }}{N}+\frac{1}{3}\frac{%
B_{2}^{\prime }}{B_{2}}-\frac{2}{3}\left( A_{1}+A_{2}\right) +\frac{2}{3}%
\frac{\left( \frac{h}{B_{2}}\right) ^{\prime }}{\frac{h}{B_{2}}}.
\end{equation*}%
Hence 
\begin{equation*}
2\left( A_{1}+A_{2}\right) =-3\frac{\left( 4q_{1}q_{3}-q_{2}^{2}\right)
^{\prime }}{\left( 4q_{1}q_{3}-q_{2}^{2}\right) }+4\frac{N^{\prime }}{N}+4%
\frac{B_{2}^{\prime }}{B_{2}}+2\frac{\left( \frac{h}{B_{2}}\right) ^{\prime }%
}{\frac{h}{B_{2}}}.
\end{equation*}%
By the definition of $A_{i}$ $\left( i=1,2\right) $ and simple integration,
we deduce that 
\begin{equation*}
\deg \left( \alpha +\beta \right) <\deg \alpha =\deg \beta 
\end{equation*}%
which is a contradiction.

\quad

\noindent \textbf{Subcase 2.} If $4q_{1}q_{3}\equiv q_{2}^{2},$ then from $%
\left( 3.29\right) $ and $\left( 3.21\right) $ we have 
\begin{equation}
\frac{B_{2}^{\prime }}{B_{2}}\frac{q_{1}}{q_{2}}-\frac{q_{1}^{\prime }}{q_{2}%
}=\frac{H}{3B_{2}N}.  \tag{3.30}
\end{equation}%
On the other hand%
\begin{equation}
\frac{q_{1}}{q_{3}}-\frac{M^{\prime }-A_{2}M}{N}=\frac{H}{3B_{2}N}. 
\tag{3.31}
\end{equation}%
\ Combining $\left( 3.30\right) $ and $\left( 3.31\right) ,$ we obtain 
\begin{equation*}
\frac{5}{4}\frac{B_{2}^{\prime }}{B_{2}}\frac{\left( \left(
A_{1}-A_{2}\right) \frac{B_{1}}{B_{2}}-\left( \frac{B_{1}}{B_{2}}\right)
^{\prime }\right) ^{\prime }}{\left( A_{1}-A_{2}\right) \frac{B_{1}}{B_{2}}%
-\left( \frac{B_{1}}{B_{2}}\right) ^{\prime }}
\end{equation*}%
\begin{equation*}
+\frac{1}{6}\left( \frac{\left( \left( A_{1}-A_{2}\right) \frac{B_{1}}{B_{2}}%
-\left( \frac{B_{1}}{B_{2}}\right) ^{\prime }\right) ^{\prime }}{\left(
A_{1}-A_{2}\right) \frac{B_{1}}{B_{2}}-\left( \frac{B_{1}}{B_{2}}\right)
^{\prime }}\right) ^{\prime }-\left( \frac{1}{2}A_{1}+A_{2}\right) \frac{%
\left( \left( A_{1}-A_{2}\right) \frac{B_{1}}{B_{2}}-\left( \frac{B_{1}}{%
B_{2}}\right) ^{\prime }\right) ^{\prime }}{\left( A_{1}-A_{2}\right) \frac{%
B_{1}}{B_{2}}-\left( \frac{B_{1}}{B_{2}}\right) ^{\prime }}
\end{equation*}%
\begin{equation*}
-\frac{5}{3}\frac{h^{\prime }}{h}\frac{\left( \left( A_{1}-A_{2}\right) 
\frac{B_{1}}{B_{2}}-\left( \frac{B_{1}}{B_{2}}\right) ^{\prime }\right)
^{\prime }}{\left( A_{1}-A_{2}\right) \frac{B_{1}}{B_{2}}-\left( \frac{B_{1}%
}{B_{2}}\right) ^{\prime }}-\frac{5}{6}\left( A_{1}+A_{2}\right) \frac{%
h^{\prime }}{h}+\frac{23}{12}\frac{B_{2}^{\prime }}{B_{2}}\frac{h^{\prime }}{%
h}-\frac{5}{4}\left( \frac{h^{\prime }}{h}\right) ^{2}
\end{equation*}%
\begin{equation*}
-\frac{1}{9}\left( A_{1}+A_{2}\right) \frac{\left( \left( A_{1}-A_{2}\right) 
\frac{B_{1}}{B_{2}}-\left( \frac{B_{1}}{B_{2}}\right) ^{\prime }\right)
^{\prime }}{\left( A_{1}-A_{2}\right) \frac{B_{1}}{B_{2}}-\left( \frac{B_{1}%
}{B_{2}}\right) ^{\prime }}+\frac{1}{9}\frac{B_{2}^{\prime }}{B_{2}}\frac{%
\left( \left( A_{1}-A_{2}\right) \frac{B_{1}}{B_{2}}-\left( \frac{B_{1}}{%
B_{2}}\right) ^{\prime }\right) ^{\prime }}{\left( A_{1}-A_{2}\right) \frac{%
B_{1}}{B_{2}}-\left( \frac{B_{1}}{B_{2}}\right) ^{\prime }}
\end{equation*}%
\begin{equation*}
+\frac{2}{9}\frac{B_{2}^{\prime }}{B_{2}}\left( A_{1}+A_{2}\right) -\frac{7}{%
9}\left( \frac{B_{2}^{\prime }}{B_{2}}\right) ^{2}-\frac{19}{36}\left( \frac{%
\left( \left( A_{1}-A_{2}\right) \frac{B_{1}}{B_{2}}-\left( \frac{B_{1}}{%
B_{2}}\right) ^{\prime }\right) ^{\prime }}{\left( A_{1}-A_{2}\right) \frac{%
B_{1}}{B_{2}}-\left( \frac{B_{1}}{B_{2}}\right) ^{\prime }}\right) ^{2}
\end{equation*}%
\begin{equation*}
-\frac{1}{6}\left( \frac{B_{2}^{\prime }}{B_{2}}\right) ^{\prime }-\frac{1}{2%
}A_{1}^{\prime }+\frac{1}{3}\left( A_{1}^{\prime }+A_{2}^{\prime }\right) +%
\frac{1}{3}\left( A_{1}+A_{2}\right) \frac{B_{2}^{\prime }}{B_{2}}+\frac{1}{2%
}A_{2}A_{1}
\end{equation*}%
\begin{equation*}
=\frac{1}{9}\left( A_{1}+A_{2}\right) ^{2}.
\end{equation*}%
Dividing both sides of the above equation by $\frac{\left(
A_{1}+A_{2}\right) ^{2}}{2}$ and since $\underset{z\rightarrow \infty }{\lim 
}\frac{R^{\prime }\left( z\right) }{R\left( z\right) }=0$ if $R$ is a
nonzero rational function, we obtain%
\begin{equation}
\left\vert \frac{A_{2}A_{1}}{\left( A_{1}+A_{2}\right) ^{2}}-\frac{2}{9}%
\right\vert \leq \frac{5}{3}\frac{\left\vert \frac{h^{\prime }}{h}%
\right\vert }{\left\vert A_{1}+A_{2}\right\vert }+\frac{23}{6}\left\vert 
\frac{B_{2}^{\prime }}{B_{2}}\right\vert \frac{\left\vert \frac{h^{\prime }}{%
h}\right\vert }{\left\vert A_{1}+A_{2}\right\vert ^{2}}+\frac{5}{2}\frac{%
\left\vert \frac{h^{\prime }}{h}\right\vert ^{2}}{\left\vert
A_{1}+A_{2}\right\vert ^{2}}+o\left( 1\right)  \tag{3.32}
\end{equation}%
On the other hand, since $\rho \left( h\right) \leq \rho \left( f\right) -1$
and by Lemma 2.3%
\begin{equation}
\left\vert \frac{h^{\prime }\left( z\right) }{h\left( z\right) }\right\vert
\leq \left\vert z\right\vert ^{\rho \left( f\right) -2+\varepsilon } 
\tag{3.33}
\end{equation}%
for all $z$ satisfying $\left\vert z\right\vert \notin E_{0}\cup \left[ 0,1%
\right] ,$ where $E_{0}\subset \left( 1,\infty \right) $ is a set of finite
logarithmic measure. By combining $\left( 3.32\right) $ and $\left(
3.33\right) ,$ we deduce 
\begin{equation*}
\underset{\underset{\left\vert z\right\vert \notin E_{0}\cup \left[ 0,1%
\right] }{z\rightarrow \infty }}{\lim }\frac{A_{2}A_{1}}{\left(
A_{1}+A_{2}\right) ^{2}}=\underset{\underset{\left\vert z\right\vert \notin
E_{0}\cup \left[ 0,1\right] }{z\rightarrow \infty }}{\lim }\frac{\alpha
^{\prime }\beta ^{\prime }}{\left( \alpha ^{\prime }+\beta ^{\prime }\right)
^{2}}=\frac{2}{9}.
\end{equation*}%
By setting $\alpha \left( z\right) =a_{m}z^{m}+\cdots +a_{0}$ and $\beta
\left( z\right) =b_{m}z^{m}+\cdots +b_{0},$ we deduce 
\begin{equation*}
\underset{\underset{\left\vert z\right\vert \notin E_{0}\cup \left[ 0,1%
\right] }{z\rightarrow \infty }}{\lim }\frac{\alpha ^{\prime }\beta ^{\prime
}}{\left( \alpha ^{\prime }+\beta ^{\prime }\right) ^{2}}=\frac{a_{m}b_{m}}{%
\left( a_{m}+b_{m}\right) ^{2}}=\frac{2}{9}
\end{equation*}%
which implies that $\frac{a_{m}}{b_{m}}=2$ or $\frac{1}{2}.$ We consider
first the case $\frac{a_{m}}{b_{m}}=\frac{1}{2},$ we get from $\left(
3.1\right) $ and $\left( 3.17\right) $ 
\begin{equation}
\varphi f^{2}+\psi f^{\prime }f-q=Ae^{\frac{1}{2}b_{m}z^{m}}  \tag{3.34}
\end{equation}%
and 
\begin{equation}
Mf^{2}+Nf^{\prime }f-q=Be^{b_{m}z^{m}},  \tag{3.35}
\end{equation}%
where $A=p_{1}e^{a_{m-1}z^{m-1}+\cdots +a_{0}}$ and $%
B=p_{2}e^{b_{m-1}z^{m-1}+\cdots +b_{0}}.$ From $\left( 3.34\right) $ and $%
\left( 3.35\right) ,$ we get%
\begin{equation*}
\varphi f^{2}+\psi f^{\prime }f=q+A\left( \frac{Mf^{2}+Nf^{\prime }f-q}{B}%
\right) ^{\frac{1}{2}}.
\end{equation*}%
Hence 
\begin{equation*}
\varphi f+\psi f^{\prime }=\frac{q}{f}+A\left( \frac{Mf^{2}+Nf^{\prime }f-q}{%
Bf^{2}}\right) ^{\frac{1}{2}}.
\end{equation*}%
Therefore 
\begin{equation*}
T\left( r,\varphi f+\psi f^{\prime }\right) =m\left( r,\varphi f+\psi
f^{\prime }\right) +S\left( r,f\right)
\end{equation*}%
\begin{equation*}
=\frac{1}{2\pi }\int_{E_{1}}\log ^{+}\left\vert \varphi \left( re^{i\theta
}\right) f\left( re^{i\theta }\right) +\psi \left( re^{i\theta }\right)
f^{\prime }\left( re^{i\theta }\right) \right\vert d\theta
\end{equation*}%
\begin{equation*}
+\frac{1}{2\pi }\int_{E_{2}}\log ^{+}\left\vert \varphi \left( re^{i\theta
}\right) f\left( re^{i\theta }\right) +\psi \left( re^{i\theta }\right)
f^{\prime }\left( re^{i\theta }\right) \right\vert d\theta +S\left(
r,f\right) ,
\end{equation*}%
where $E_{1}=\left\{ \theta :\left\vert f\left( re^{i\theta }\right)
\right\vert \leq 1\right\} $ and $E_{2}=\left\{ \theta :\left\vert f\left(
re^{i\theta }\right) \right\vert >1\right\} .$ Now%
\begin{equation*}
\frac{1}{2\pi }\int_{E_{1}}\log ^{+}\left\vert \varphi \left( re^{i\theta
}\right) f\left( re^{i\theta }\right) +\psi \left( re^{i\theta }\right)
f^{\prime }\left( re^{i\theta }\right) \right\vert d\theta
\end{equation*}%
\begin{equation*}
\leq \frac{1}{2\pi }\int_{E_{1}}\log ^{+}\left\vert f^{\prime }\left(
re^{i\theta }\right) \right\vert d\theta +S\left( r,f\right)
\end{equation*}%
\begin{equation*}
\leq \frac{1}{2\pi }\int_{E_{1}}\log ^{+}\left\vert \frac{f^{\prime }\left(
re^{i\theta }\right) }{f\left( re^{i\theta }\right) }\right\vert d\theta
+S\left( r,f\right) =S\left( r,f\right) .
\end{equation*}%
On the other hand 
\begin{equation*}
\frac{1}{2\pi }\int_{E_{2}}\log ^{+}\left\vert \varphi \left( re^{i\theta
}\right) f\left( re^{i\theta }\right) +\psi \left( re^{i\theta }\right)
f^{\prime }\left( re^{i\theta }\right) \right\vert d\theta
\end{equation*}%
\begin{equation*}
=\frac{1}{2\pi }\int_{E_{2}}\log ^{+}\left\vert \frac{q\left( re^{i\theta
}\right) }{f\left( re^{i\theta }\right) }\right\vert d\theta
\end{equation*}%
\begin{equation*}
+\frac{1}{4\pi }\int_{E_{2}}\log ^{+}\left\vert \frac{M\left( re^{i\theta
}\right) }{B\left( re^{i\theta }\right) }+\frac{N\left( re^{i\theta }\right) 
}{B\left( re^{i\theta }\right) }\frac{f^{\prime }\left( re^{i\theta }\right) 
}{f\left( re^{i\theta }\right) }-\frac{q\left( re^{i\theta }\right) }{%
f^{2}\left( re^{i\theta }\right) }\right\vert d\theta +S\left( r,f\right)
=S\left( r,f\right) .
\end{equation*}%
Hence 
\begin{equation*}
T\left( r,f_{c_{1}}\right) =T\left( r,\varphi f+\psi f^{\prime }\right)
=S\left( r,f\right)
\end{equation*}%
which is a contradiction. If $\frac{a_{m}}{b_{m}}=2,$ then by the same
argument we have 
\begin{equation*}
Mf^{2}+Nf^{\prime }f=q+B\left( \frac{\varphi f^{2}+\psi f^{\prime }f-q}{A}%
\right) ^{\frac{1}{2}}
\end{equation*}%
which implies the contradiction 
\begin{equation*}
T\left( r,f_{c_{2}}\right) =T\left( r,Mf+Nf^{\prime }\right) =S\left(
r,f\right) .
\end{equation*}%
\textbf{Case 2.} $B_{2}f_{c_{1}}-B_{1}f_{c_{2}}\equiv 0,$ by using the same
arguments as in the proof of $\left( 3.14\right) ,$ we obtain that%
\begin{equation*}
\frac{A_{1}B_{1}}{B_{2}}-\left( \frac{B_{1}}{B_{2}}\right) ^{\prime }-\frac{%
B_{1}A_{2}}{B_{2}}\equiv 0
\end{equation*}%
which leads to 
\begin{equation}
\frac{p_{1}}{p_{2}}e^{\alpha -\beta }=k\frac{B_{1}}{B_{2}}=k\frac{f_{c_{1}}}{%
f_{c_{2}}},  \tag{3.36}
\end{equation}%
where $k$ is a nonzero complex constant. By this $\left( 3.1\right) $ and $%
\left( 3.2\right) ,$ we have%
\begin{equation}
\left( 1-c\right) ff_{c_{1}}f_{c_{2}}=q\left( f_{c_{2}}-kf_{c_{1}}\right) . 
\tag{3.37}
\end{equation}%
If $k\neq 1,$ then by applying Clunie lemma to $\left( 3.37\right) ,$ we
deduce the contradiction $T\left( r,f_{c_{i}}\right) =S\left( r,f\right) .$
Hence, $k=1$ and from the equation $\left( 3.36\right) ,$ we conclude that $%
f_{c_{1}}\equiv f_{c_{2}}$ which exclude the hypothesis of our theorem. This
shows that at least one of $f\left( z\right) f\left( z+c_{1}\right) -q\left(
z\right) $ and $f\left( z\right) f\left( z+c_{2}\right) -q\left( z\right) $
has infinitely many zeros.

\quad

\noindent \textbf{Acknowledgements.} The authors are grateful to the
anonymous referee for his/her valuable comments and suggestions which lead
to the improvement of this paper.

\begin{center}
{\Large References}
\end{center}

\noindent $\left[ 1\right] $ W. Bergweiler, \textit{On the product of a
meromorphic function and its derivative}, Bull. Hong Kong Math. Soc. 1
(1997), 97--101.

\noindent $\left[ 2\right] $ W. Bergweiler and A. Eremenko, \textit{On the
singularities of the inverse to a meromorphic function of finite order},
Rev. Mat. Iberoamericana 11 (1995), no. 2, 355--373.

\noindent $\left[ 3\right] $ H. H. Chen and M. L. Fang, \textit{The value
distribution of }$f^{n}f^{\prime }$, China Ser. A 38 (1995), no. 7, 789--798.

\noindent $\left[ 4\right] \ $Z. X. Chen, \textit{Complex differences and
difference equations}, Mathematics Monograph series 29, Science Press,
Beijing, 2014.

\noindent $\left[ 5\right] \ $Z. X. Chen, \textit{On value distribution of
difference polynomials of meromorphic functions}. Abstr. Appl. Anal. 2011,
Art. ID 239853, 9 pp.

\noindent $\left[ 6\right] \ $Y. M. Chiang, S. J. Feng, \textit{On the
Nevanlinna characteristic of }$f\left( z+\eta \right) $ \textit{and
difference equations in the complex plane, }Ramanujan J. 16 (2008), no. 1,
105--129.

\noindent $\left[ 7\right] \ $G. G. Gundersen, \textit{Estimates for the
logarithmic derivative of a meromorphic function, plus similar estimates},
J. London Math. Soc. (2) 37 (1988), no. 1, 88--104.

\noindent $\left[ 8\right] \ $R. G. Halburd, R. J. Korhonen, \textit{%
Difference analogue of the lemma on the logarithmic derivative with
applications to difference equations, }J. Math. Anal. Appl. 314 (2006)%
\textit{, }no. 2, 477--487.

\noindent $\left[ 9\right] \ $R. G. Halburd, R. J. Korhonen, \textit{%
Nevanlinna theory for the difference operator}, Ann. Acad. Sci. Fenn. Math.
31 (2006), no. 2, 463--478.

\noindent $\left[ 10\right] \ $W. K. Hayman, \textit{Meromorphic functions},
Oxford Mathematical Monographs Clarendon Press, Oxford 1964.

\noindent $\left[ 11\right] \ $W. K. Hayman, \textit{Picard values of
meromorphic functions and their derivatives}, Ann. of Math. (2) 70 (1959),
9--42.

\noindent $\left[ 12\right] $ W. K. Hayman, \textit{Research problems in
function theory}, The Athlone Press University of London, London, 1967.

\noindent $\left[ 13\right] \ $I. Laine, \textit{Nevanlinna theory and
complex differential equations}, de Gruyter Studies in Mathematics, 15.
Walter de Gruyter \& Co., Berlin, 1993.

\noindent $\left[ 14\right] \ $I. Laine and C. C. Yang, \textit{Value
distribution of difference polynomials}, Proc. Japan Acad. Ser. A Math. Sci.
83 (2007), no. 8, 148--151.

\noindent $\left[ 15\right] \ $K. Liu and L. Z. Yang, \textit{Value
distribution of the difference operator}, Arch. Math. (Basel) 92 (2009), no.
3, 270--278.

\noindent $\left[ 16\right] $ N. Li and L. Yang, \textit{Value distribution
of certain type of difference polynomials}, Abstr. Appl. Anal. 2014, Art. ID
278786, 6 pp.

\noindent $\left[ 17\right] \ $W. L\"{u}, N. Liu, C. Yang and C. Zhuo, 
\textit{Notes on value distributions of }$ff^{\left( k\right) }-b,$ Kodai
Math. J. 39 (2016), no. 3, 500--509.

\noindent $\left[ 18\right] $ E. Mues, \textit{\"{U}ber ein Problem von
Hayman}, Math. Z. 164 (1979), no. 3, 239--259. (in German).

\end{document}